\documentclass[12pt]{article}

\usepackage{amsbsy}
\usepackage{amsmath}
\usepackage{amssymb}
\usepackage{graphicx}

\newcommand{\eh}{\hspace{.06in}}
\newcommand{\ih}{\'\i}
\newcommand{\C}{\mathbb{C}}
\newcommand{\E}{\mathbb{E}}
\newcommand{\R}{\mathbb{R}}
\newcommand{\B}{\mathcal{B}}
\newcommand{\mC}{\mathcal{C}}
\newcommand{\D}{\mathcal{D}}
\newcommand{\K}{\mathcal{K}}
\newcommand{\mP}{\mathcal{P}}
\newcommand{\sC}{\mathsf{C}}
\newcommand{\sP}{\mathsf{P}}
\newcommand{\sS}{\mathsf{S}}
\newcommand{\oM}{\overline{M}}
\newcommand{\CL}{{\rm CL}_{\frac{\pi}{2}}}
\newcommand{\deh}{\partial}
\newcommand{\ds}{\displaystyle}
\newcommand{\af}{\alpha}
\newcommand{\eps}{\varepsilon}
\newcommand{\ld}{\lambda}
\newcommand{\ta}{\theta}
\newcommand{\Om}{\Omega}
\newcommand{\BE}{\begin{equation}}
\newcommand{\EE}{\end{equation}}
\newcommand{\fns}{\footnotesize}
\newcommand{\Lim}[1]{\lower5.5pt\hbox{${{\ds\lim}\atop^{#1}}$} \ }
 
\begin{document}
\centerline{\large{\bf A LIMIT-METHOD FOR SOLVING}}
\centerline{\large{\bf PERIOD PROBLEMS ON MINIMAL SURFACES}}
\bigskip

\centerline{K{\fns ELLY} L{\fns \"UBECK} \ \& \ V{\fns AL\'ERIO} R{\fns AMOS} B{\fns ATISTA}}
\bigskip

\begin{abstract}
We introduce a new technique to solve period problems on minimal surfaces called ``limit-method''. If a family of surfaces has Weierstra\ss  data converging to the data of a known example, and this presents a transversal solution of periods, then the original family contains a sub-family with closed periods.
\end{abstract}
\ \\
\centerline{\bf 1. Introduction}
\\

The study of minimal surfaces was first motivated by their physical properties. Given a simple closed curve $\sC$ in $\R^3$, the surface $\sS$ with $\sC=\deh\sS$ and least area is also tension-minimiser. More generally, if $\sS$ separates two homogeneous means, each under uniform pressures $\sP_1$ and $\sP_2$, tension-minimising is equivalent to constant $H$, which denotes the mean curvature of $\sS$. In fact, $H$ is then proportional to $\sP_1-\sP_2$ (see \cite{BCE} for details).

In 1883, Enneper and Weierstra\ss  showed that every minimal surface is locally parametrised by $F:\overline{\Om}\to\R^3$, for open $\Om\subset\C$ and $F$ determined by a meromorphic function $g$ and a holomorphic differential $dh$ on $\Om$. One proves that $g$ is the stereographic projection of the Gau\ss  map $N:\Om\to S^2$. Moreover, if $F_\ta$ is the related to map obtained from $e^{i\ta}dh$, $\ta\in[0,2\pi)$, then all $F_\ta:\overline{\Om}\to\R^3$ give isometric minimal immersions. The case $\ta=\pi/2$ is denoted $F^*$ or $\sS^*$ and called {\it conjugate}. 

If $F(\deh\Om)=\sC$ is polygonal with bijective orthogonal projection $\mP$ onto a plane, then $F$ is the graph of a function $f:\mP\cup Int\mP\to\R$. In this case, since $f|_\mP$ is explicit, a numeric PDE-solution describes the surface by computer graph. From the isometry and the coinciding Gau\ss  maps of $F$ and $F^*$, one also gets numerical pictures of $\sS^*$. Such a procedure is called ``Conjugate Plateau Construction''.

In some cases, $\sS$ is contained in a complete surface $M$ with $H\equiv 0$, still called ``minimal'' although $\deh M=\emptyset$. For instance, the Schwarz Reflection Principle applies to polygonal $\sC$, and so $M$ is given by $X:R\to\R^3$, for a Riemann surface $R$ with local chart $\psi:\Om\to R$, and $F=X\circ\psi$. This allows $(g,dh)$ to be globally defined on $R$, and if $R$ is compact, one of its algebraic equations can give explicit formulae for $g$ and $dh$. This is the second case where computer graphs become possible, even to give a notion of the {\it whole} $M$ itself. This procedure is called ``Weierstra\ss  Data Construction''. 

To date, one has not found other explicit constructions. For surfaces called {\it algebraic}, namely complete and with finite total curvature in a flat space, implicit constructions were introduced by Traizet and Kapouleas, but it either lacks $(g,dh)$ or the underlying $R$ (see \cite{K} and \cite{T1}). The greatest difficulty with the ``Weierstra\ss  Construction'' for algebraic surfaces are the so-called {\it period problems}. In general, they are a system of equations involving elliptic integrals with several interdependent parameters. If ever solvable, it is usually with extreme difficulties. 

Other constructions could be called {\it almost explicit}, where it only lacks the possibility for a theoretical refinement of the parameters domain. This is the case of \cite{T2}, \cite{T3} and \cite{W2}, in which one knows the domain to be in a punched neighbourhood of $0\in\R^n$, however impossible to describe, where $n$ is the number of parameters. In those works, one applies the {\it implicit function theorem}, which requires the partial derivatives of elliptic integrals computed at the origin.  

We classify the method presented herein as ``almost explicit''. The difference is that practically no computation is needed to solve periods, but the sought after examples must converge to a limit-surface, of which the periods have {\it transversal} solution. Here this term just means that the difference of two continuous functions changes sign, disconsidering whether they have non-coinciding tangent hyperplanes along the crossing, as required in the classical sense.

The method not only solves periods, but also helps to prove embeddedness. Since one gets a continuous family of period-free surfaces, an embedded limit-member from the family may be used for this purpose (see \cite{HK} or \cite{W1} as a reference). Some examples with transversal solution of periods are found in \cite{CHM} or \cite{HK}. Neither Costa's nor Chen-Gackstatter surfaces fit this requirement (see \cite{C}, \cite{CG}, \cite{HK} and \cite{Ka}).      

In order to illustrate our method, we shall construct examples labelled $\CL$. These were inspired in the surfaces L$_\frac{\pi}{2}$ and C$_\frac{\pi}{2}$ from \cite{V1} and \cite{V2}. In Section 3 one computes Weierstra\ss  data by Karcher's {\it reverse construction}, a powerful method described in \cite{Ka}. In Section 4 one describes the analytic equations for the periods, which are finally solved in Section 5 by our limit-method. Section 6 is devoted to the embeddedness proof of the fundamental piece $P$, which generates $\CL$ by isometries in $\R^3$. 

This work refers to part of the first author's doctoral thesis \cite{L}, which was supported by CAPES - Coordena\c c\~ao de Aperfei\c coamento de Pessoal de N\ih vel Superior.
\\
\\
\centerline{\bf 2. Preliminaries}
\\

In this section we state some basic definitions and theorems. Throughout this work, surfaces are considered connected and regular. Details can be found in \cite{Ka}, \cite{LM}, \cite{N} and \cite{O}.

{\bf Theorem 2.1.} \it Let $X:R\to\E$ be a complete isometric immersion of a Riemannian surface $R$ into a three-dimensional complete flat space $\E$. If $X$ is minimal and the total Gaussian curvature $\int_R K dA$ is finite, then $R$ is biholomorphic to a compact Riemann surface $\overline{R}$ punched at a finite number of points.\rm

{\bf Theorem 2.2.} (Weierstra\ss  representation). \it Let $R$ be a Riemann surface, $g$ and $dh$ meromorphic function and 1-differential form on $R$, such that the zeros of $dh$ coincide with the poles and zeros of $g$. Suppose that $X:R\to\E$, given by
\[
   X(p):=Re\int^p(\phi_1,\phi_2,\phi_3),\eh\eh where\eh\eh
   (\phi_1,\phi_2,\phi_3):=\frac{1}{2}(g^{-1}-g,ig^{-1}+ig,2)dh,
\]
is well-defined. Then $X$ is a conformal minimal immersion. Conversely, every conformal minimal immersion $X:R\to\E$ can be expressed as above for some meromorphic function $g$ and 1-form $dh$.\rm

{\bf Definition 2.1.} The pair $(g,dh)$ is the \it Weierstra\ss  data \rm and $\phi_{1,2,3}$ are the \it Weierstra\ss  forms \rm on $R$ of the minimal immersion $X:R\to X(R)\subset\E$.

{\bf Theorem 2.3.} \it Under the hypotheses of Theorems 2.1 and 2.2, the Weierstra\ss  data $(g,dh)$ extend meromorphically on $\overline{R}$.\rm

The function $g$ is the stereographic projection of the Gau\ss \ map $N:R\to S^2$ of the minimal immersion $X$. It is a branched covering of $\hat\C$ and $\int_SKdA=-4\pi$deg$(g)$.
\\
\\
\centerline{\bf 3. The Weierstra\ss  data of $\CL$}
\\

Let us consider a topological surface which could admit a Riemannian structure and be isometrically immersed in $\R^3$ with the following properties: 1) the immersed surface is minimal and doubly periodic; 2) it is spanned by reflection in a vertical plane, together with a horizontal translation group, both applied to a fundamental piece $P$, where $P$ is a surface with boundary and two catenoidal ends; 3) $P$ has a symmetry group generated by $180^\circ$-rotations about two segments crossing orthogonally at their middle point $S$; 4) $\deh P$ consists of two congruent curves alternating with two segments, and they project onto a rectangle $Q$. Figure 1 illustrates the sought after surface. 

\begin{figure}[ht]
\includegraphics[scale=0.55]{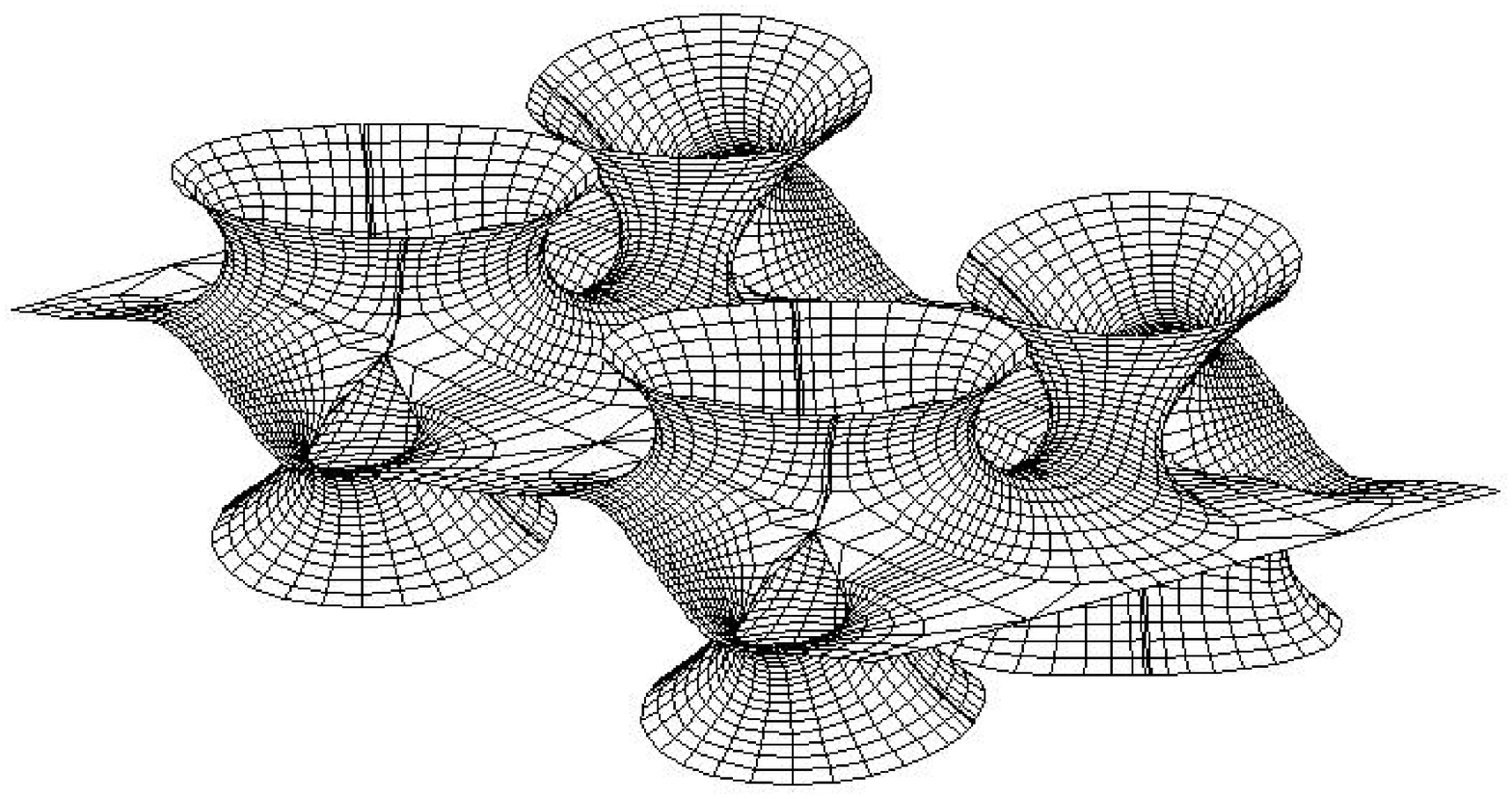}
\caption{The surfaces $\CL$.}
\label{f1}
\end{figure}

One can interpret $P$ as the Costa surface with its catenoidal ends kept, but the planar end replaced by $\deh P$. Point $S$ represents the ``Costa-saddle'', where the lines meet. Of course, the spanned doubly periodic surface will have self-intersections, but they will occur only at the catenoidal ends. Let us consider $S$ as the origin of $\R^3$ and the segments of $P\setminus\deh P$ contained in $Ox_1$, $Ox_2$. Therefore, $Ox_3$ becomes the axis for both ``top'' and ``bottom'' catenoidal ends. Moreover, we consider that $x_3=x_1\wedge x_2$ and $Ox_1$ intersects only the straight segments of $\deh P$. By projecting $\deh P$ orthogonally onto $x_3=0$, one sees that $Q$ must be orthogonal to $Ox_{1,2}$. Otherwise, it would occur more self-intersections than just at the ends. 

\begin{figure}[ht]
\centering
\includegraphics[scale=1.2]{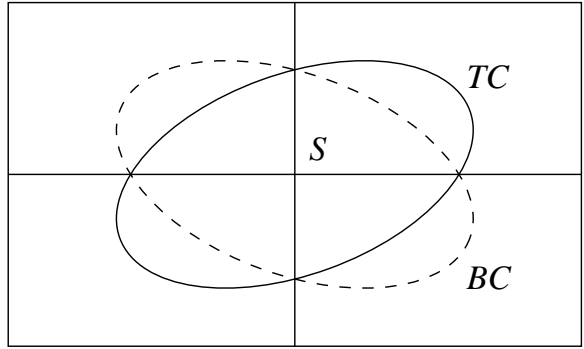}
\caption{The rectangle $Q$.}
\end{figure}

From Figure 2 we notice that $Q$ and the ends can assume different ratios and {\it logarithmic growths}, respectively (see \cite{HK} for a definition). Therefore, the $\CL$ examples make a two-parameter family of surfaces. Let us take a generic member from this family. Consider its quotient by the translation group followed by a compactification of the four ends. We then get a topological compact surface denoted $\oM$ (see Figure 3a). From this picture one easily sees that $\oM$ has genus 3. The stretch $A-L-F-D-A$ generates the surface lines, whereas $A-E-A$ and $F-N-F$ generates the reflectional symmetry curves.

Let us now consider $\rho:\oM\to\oM$ as the involution corresponding to $180^{\circ}$ rotation about $Ox_3$. It fixes $L,D,TC,S,BC,TC',S'$ and $BC'$, but interchanges $A,F$ and $E,N$. Hence, the Euler-Poincar\'e formula gives
\[
   \chi(\oM/\rho)=\frac{\chi(\oM)}{2}+\frac{8}{2}=2.
\]
Namely, $\oM/\rho$ is topologically $S^2$, which admits the unique conformal structure $\hat{\C}=\C\cup\{\infty\}$ by Koebe's theorem. Since $\rho$ is a branched covering, it induces a conformal structure on $\oM$. This determines a meromorphic map $z:\oM\to\hat{\C}$, and up to a M\"obius transformation we choose $z(S)=0$, $z(A)=1$ and $z(S')=\infty$. Rotation of $180^{\circ}$ about $A-L$ fixes $0$, $1$ and $\infty$. Therefore, the induced involution in $\hat{\C}$ is a reflection in the real line. Consequently, $z(A-L)\subset\R$. However, neither $S$ nor $S'$ belong to $A-L$. This means that $z(A-L)=[1,\tilde{x}]$ for some $\tilde{x}>1$ or $z(A-L)=[x,1]$ for a certain $x<1$. Since $S-L$ is fixed by $180^{\circ}$-rotation about itself, and this fixes $S'$, then $z(S-L)\subset\R$. This implies that $z(S)=0<z(L)<z(A)=1$ and so $z(L)=x$ with $0<x<1$.

That rotation of $S-L$ about itself also fixes $L$ and $D$, namely $0$, $x$, $\infty$ and $z(D)$. So, it induces $z\to\bar{z}$ in $\hat{\C}$. But it interchanges pairs $(A,F)$, $(E,N)$, $(TC,BC)$ and $(TC',BC')$, while $z(E)=z(N)$ by $\rho$. Therefore,
\BE
   \bar{z}(N)=z(E),\eh{\rm which\eh implies}\eh z(E)=z(N)\in\R,
\EE
and
\BE
   z(TC)=\bar{z}(BC),\eh\bar{z}(TC')=z(BC').
\EE

\begin{figure}[!htb]
\begin{minipage}[b]{0.40\linewidth}
\includegraphics[width=\linewidth]{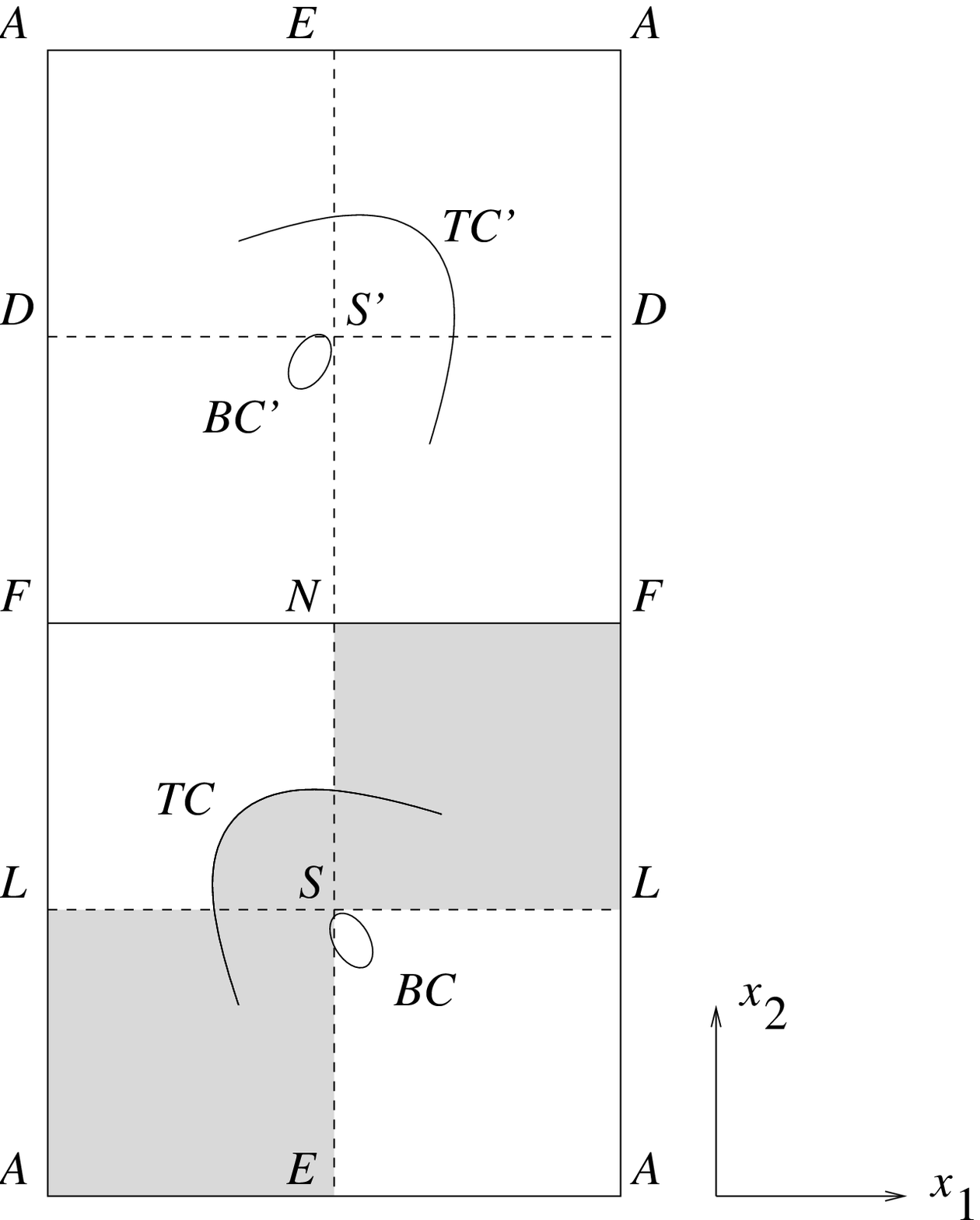}
\end{minipage}\hfill
\begin{minipage}[b]{0.50\linewidth}
\includegraphics[width=\linewidth]{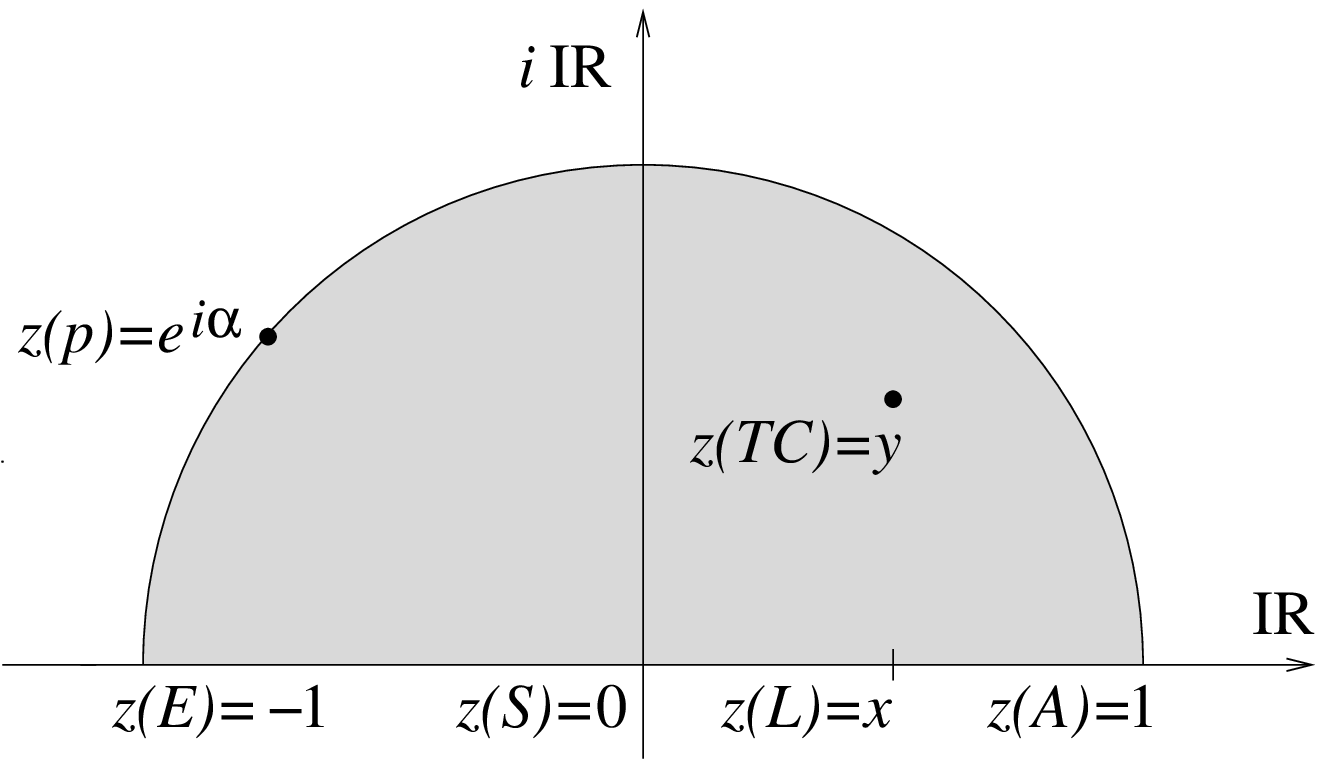}
\end{minipage}
\caption{(a) The surface $\oM$ with important points; (b) the map $z:\oM\to\hat{\C}$.}
\end{figure}

Reflection in $A-E$ fixes $A,E,F,N$ and interchanges pairs $(TC,TC')$, $(BC,BC')$, $(S,S')$ and $(L,D)$. This means, it fixes $1$ and interchanges $0,\infty$. Therefore, the induced involution in $\hat{\C}$ is given by $z\to1/\overline{z}$, whence $z(D)=1/\bar{z}(L)=1/x$ and
\BE
   z(TC)=1/\bar{z}(TC'),\eh z(BC)=1/{z}(BC').
\EE

Observe that the fixed points of this involution belong to $S^{1}$. Therefore
\BE
   z(E)=z(N)\in S^1. 
\EE

From (1), (4) and the fact that $z(A)=1$ we get $z(E)=z(N)=-1$. Since rotation about $A-L$ does not fix $TC$, then $z(TC)=:y\in\C\setminus\R$. From (2) and (3) we get $z(BC)=\bar{y}$, $z(TC')=1/\bar{y}$ and $z(BC')=1/y$. The curves $A-E-A$ and $F-N-F$ are in planes parallel to $x_2=0$ and have highest and lowest points. Let $p,p'$ and $q,q'$ be such points on $A-E-A$ and $F-N-F$, respectively. Notice that $\measuredangle(A-E,A-L)=\pi/2$ and $A$ is not fixed by $\rho$, thus  $\measuredangle(z(A-E),z(A-L))=\pi/2$. Since reflection in $A-E$ is given by $z\to 1/\bar{z}$, then $z(A-E)\subset S^{1}$. Without loss of generality we can take $Im\{y\}>0$, which implies $|y|<1$ and $z(p)=e^{i\af}$ for some $\af\in(0,\pi)$. Notice that $\rho$ interchanges $p$ with $q$, thus $z(q)=e^{i\af}$. Rotation about $S-E$ gives $z(p')=\bar{z}(p)$ and $z(q^\prime)=\bar{z}(q)$, hence $z(p')=z(q')=e^{-i\af}$. Figure 3b illustrates the map $z:\oM\to\C$. 

Since $\rho:\oM\to\oM$ is the hyperelliptic involution, it is easy to write down an algebraic equation for $\oM$:
\BE
   w^2=\frac{Z^{2} -2Re\{Y\}Z +|Y|^{2}}{(Z-X)(Z-2\cos\af)^{2}},
\EE
where $X=x+1/x$, $Y=y+1/y$ and $Z=z+1/z$ with $x\in(0,1),|y|<1,Im\{y\}>0$ and $\af\in(0,\pi)$. The values $0^{\pm 1},x^{\pm 1},y^{\pm 1}$ and $\bar{y}\,^{\pm 1}$ give exactly all branch points of $z$, each of order $1$. From Riemann-Hurwitz formula, $8/2-2+1=3$, which agrees with the expected genus of $\oM$. Now suppose that $M$ is a complete minimal immersion of $\oM\setminus\{TC,BC,TC',BC'\}$ in $\R^3$. For the Weierstra\ss  data $(g,dh)$, since $g$ is the stereographic projection of the unitary normal on $M$, based on Figures 1 and 3a we settle
\BE
   g=icw,\text{ \ where \ } c>0.
\EE

From the sought after ends and regular points of $M$, $g$ determines all zeros and poles of $dh$, including their branch order. Therefore, we know that $S,S',D,L,p$, $p',q$ and $q'$ are exactly the points where $dh$ vanishes, whereas it takes $\infty$ precisely at $TC,BC,TC'$ and $BC'$. By comparing $dh$ with $z$ and $dz$, we settle
\BE
   dh=-\frac{i(Z-2\cos\af)dz/z}{Z^2-2Re\{Y\}Z +|Y|^2}.
\EE

Table (8) summarises the involutions of $\oM$ and symmetries of $M$:
\BE
\begin{tabular}{|c|c|c|c|c|}\hline 
Symmetry & Involution                    & $z$          & $w$   & $dh$  \\\hline
$E-S$    & $(w,z)\to(-\bar{w},\bar{z})$  & $-1<\cdot<0$ & $i\R$ & $i\R$ \\\hline
$S-L$    & $(w,z)\to(\bar{w},\bar{z})$   & $0<\cdot<x$  & $\R$  & $i\R$ \\\hline
$L-A$    & $(w,z)\to(-\bar{w},\bar{z})$  & $x<\cdot<1$  & $i\R$ & $i\R$ \\\hline
$A-E$    & $(w,z)\to(-\bar{w},1/\bar{z})$& $|\cdot|=1$  & $i\R$ & $\R$  \\\hline
\end{tabular}
\EE

According to this table, $dh\cdot dg/g$ is real only on $A-E$, and purely imaginary otherwise. Hence, the sought after surfaces really have all the expected symmetry curves and lines. 
\\
\\
\centerline{\bf 4. Period Analysis}
\\

Figure 4 represents the fundamental piece $P$ of $M$ with open periods. There we consider the curve $\gamma:=S\to L\to A\to E\to S$.

\begin{figure}[ht]
\centering
\includegraphics[scale=0.6]{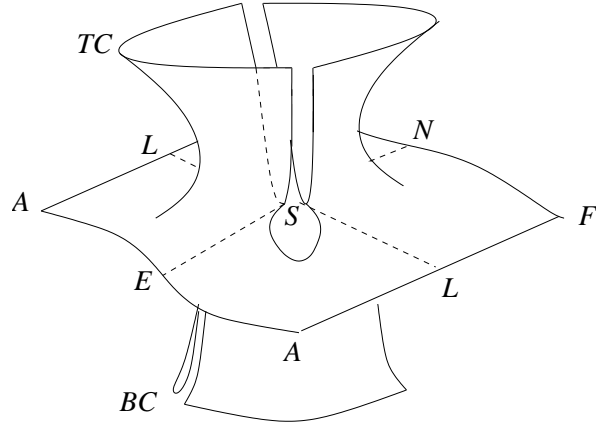}
\caption{The fundamental piece with open periods.}
\end{figure}

For every closed curve in $\oM\setminus\{TC,BC,TC',BC'\}$, one analyses the period vector $Re\oint\phi_{1,2,3}$. For instance, if such a curve is homotopic to $TC$ in $\oM$, then it is homotopically $\rho(\gamma)\cup\gamma$. Hence, one {\it must} prove that $Re\int_{\gamma}\phi_{1,2,3}=0$. The same conclusion holds for $BC$, $TC'$ and $BC'$, after applying involutions of $\oM$. For $A\to L\to F$, the period is zero on $Ox_{1,3}$ but non-zero on $Ox_2$, since it is taken to a line parallel to that axis, under the minimal immersion. Now, $A\to E\to A$ is in a plane parallel to $Ox_1x_3$, whence its period vanishes on $Ox_2$. If the period is zero on $\gamma$, we shall have $Re\int_{A\to E}\phi_1=Re\int_{L\to S}\phi_1\not=0$, for the latter is a straight segment in $M$. Finally, it is {\it enough} to prove that the period is zero on $\gamma$. Namely, the constant $c=c(x,y,\af)$ and $\af=\af(y)$ must satisfy the following three conditions:
\BE
   \cos\af=\frac{\int_0^\pi\frac{\cos tdt}{4\cos^2t-4Re\{Y\}\cos t+|Y|^2}}{\int_0^\pi\frac{dt}{4\cos^2t-4Re\{Y\}\cos t+|Y|^2}};\eh c^2=c_1:=\frac{I_1+I_2}{-I_3+I_4};\eh c^2=c_2:=\frac{I_5-I_6}{I_7-I_8};
\EE
\ \\
where $I_1=\int_0^x|dh/w|$, $I_3=\int_0^x|wdh|$, $I_5=\int_{-1}^0|dh/w|$, $I_6=\int_x^1|dh/w|$, $I_7=\int_x^1|wdh|$ and $I_8=\int_{-1}^0|wdh|$, for $z(t)=t$ varying in real intervals. Regarding $I_{2,4}$, $z(t)=e^{it}$ and $0\le t\le\pi$, so that $I_2=\int_0^\pi|dh/w|$ and $I_4=\int_0^\pi|wdh|$. Notice that $I_{2,4}$ remain invariant if one chooses $z(t)=e^{-it}$ instead. This fact will be used in the next section.

From (5) and (7), one clearly recognises the complexity of Equations (9). If we tried the {\it intermediate value theorem}, many cares would be necessary. For instance, one needs to survey the $(x,y)$-region where both denominators of $c_{1,2}$ do not vanish. Afterwards, positiveness must hold in a certain connected subregion, of which the boundary has points where $c_1-c_2$ changes sign close by. For Equations (9), the authors realised that these procedures were far too laborious and fruitless. In the next section, we apply the {\it limit-method} to solve (9), practically without computations.
\\
\\
\centerline{\bf 5. Application of the Limit-Method}
\\

In order to apply the method explained in the Introduction, we shall first analyse the Weierstra\ss  data of $M=M_{x,y,\af}$. If one considers the function $z$ as a variable in the complex plane, one of the limits for $x\to 0$ will give the Weierstra\ss  data of the so-called $M_{\rm L_b}$-surfaces (see \cite{V3}, p482). For this latter we have a solution of periods given by the transversal crossing of two graphs (see \cite{V3}, pp486-9). Roughly saying, these are graphs of $xc_{1,2}$ at $x=0$, which means that the functions $c_{1,2}$ coincide for small $x>0$, at certain values of $y$ and $\af$, which will depend on $x$. Moreover, the crossing happens at positive values of $xc_{1,2}|_{x=0}$, so that $c^2$ is positive.

Take any $\ld>1$ and $\rho\in(-\pi/2,0]$. If $x$ is sufficiently close to zero, the choice
\BE
   \bar{y}(x):=\frac{xe^{i\rho}}{e^{i\rho}+i\ld}
\EE
guarantees that $\bar{y}(x)$ belongs to the interior of $\D^-:=\{z\in\C:|z|\le 1\le 1-Im\{z\}\}$.
\begin{figure}[ht]
\centering
\includegraphics[scale=0.6]{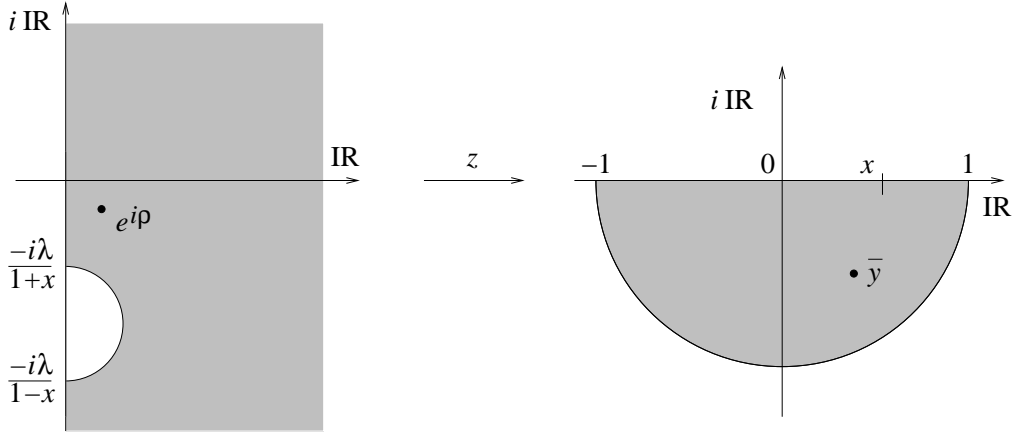}
\caption{The map $z(\zeta)$.}
\end{figure}

Consider $C:=\sqrt{x\ld}\cdot c$ and $K$ a compact set in the region $\B:=\{\zeta\in\hat{\C}:Re\{\zeta\}\geq 0\}\setminus\{-\ld i,e^{i\rho}\}$. Soon we shall prove that $C^2$ has a finite positive limit, no matter which equation one takes from (9). Assume this for the moment and also take $C_j:=x\ld\cdot c_j$, $j=1,2$. For $x$ sufficiently close to zero, the sector $\{\zeta\in\C:|\zeta+\frac{i\lambda}{1-x^{2}}|<\frac{\lambda x}{1-x^{2}}\}$ will be disjoint from $K$. Let $z:K\to\D^-$ be the map $z(\zeta)=\frac{x\zeta}{\zeta+i\ld}$ and define $G(\zeta):=g\circ z(\zeta)$, $dH:=\frac{\ld}{x}z^*dh$. Now fix any $\zeta\in K$. From (5), (6) and (7) we have 
\BE
   \lim_{x\to0} G^2(\zeta)=\lim_{x\to 0}\frac{-C^2w^2(\zeta)}{x\ld}=
   \frac{iC^2\zeta(\zeta-e^{i\rho})(\zeta+e^{-i\rho})}{(\zeta+i\ld)^2};
\EE
\BE
   \lim_{x\to 0} dH(\dot{\zeta})=\lim_{x\to 0}\frac{\ld dh(\dot{\zeta})}{x}=
   \frac{d\zeta}{(\zeta-e^{i\rho})(\zeta+e^{-i\rho})}.
\EE

One readily recognises (11) and (12) as the Weierstra\ss  data of $M_{\rm L_b}$, up to $90^\circ$-rotation about the vertical axis (see \cite{V3}, pp483-4). From this and (9), it easily follows that $\sqrt{\frac{\ld^3}{x}}I_{1,5,6}$ and $\sqrt{\frac{\ld}{x^3}}I_{3,7,8}$ will converge to the integrals in (13) and (22) of \cite{V3}, pp486-8. Now, {\it uniform} convergence can be established by the following argument. 

Notice that $I_5-I_6=\int_\sigma dh/w$, where $\sigma$ is homotopic to $\deh\D^-\setminus(0,x)$ in $\D^-\setminus\{\bar{y}\}$. Similarly, $I_7-I_8=\int_\sigma wdh$. A suitable choice of $K$ will guarantee that $(\{0\}\times i[0,+\infty])\cup\{z^*\sigma\}\subset K$. From \cite{V3}, pp485-8, one sees that the corresponding limit-integrals share the same property. Convergence for $x$ approaching zero will then be uniform, because the integrals are well-defined on paths contained in $K$. 

It remains to analyse $I_{2,4}$. In the sequel we prove that, for a curve $\zeta(t)$ homotopic to $\beta(t):=\frac{-i\ld}{1-x^2}(1+xe^{it})$, $t\in[0,\pi]$, one has $\Lim{x\to 0}{\sqrt{\frac{\ld^3}{x}}}I_2=0$ and  $\Lim{x\to 0}{\sqrt{\frac{\ld}{x^3}}}I_{4}=\frac{\ld \pi}{f(\ld)}$, where $f(\ld)=\sqrt{\ld(1+\ld^2+2\ld\sin\rho)}$. From this fact and the uniform convergence of the Weierstra\ss  data from $M_{x,y,\af}$ to $M_{\rm L_b}$ in $K$, it follows that $C_1$ and $C_2$ coincide at $x=0$ with the constants (13) and (22) of \cite{V3}, pp486-8. Namely, $C_{1,2}$ are transversal in a neighbourhood of $x=0$. Therefore, $C_1=C_2$ (and consequently $c_1=c_2$) for $x$ close to zero, $\ld=\ld(x)$ with $\ld(0)=\ld_\rho$, and $\bar{y}(x)=xe^{i\rho}/(e^{i\rho}+i\ld)$.

On $\beta(t)$, a careful computation shows that $\Lim{x\to 0}{z(t)}=e^{-it}$. Therefore,
\BE
   \lim_{x\to 0}\sqrt{x}\cdot|w|\biggl|_{\beta(t)}=\frac{|i\ld+e^{i\rho}|}{2|\cos t-\cos\af|}
\EE 
and
\BE
  \lim_{x\to 0}\frac{|dh|}{x^2}\biggl|_{\beta(t)}=\frac{2|\cos t-\cos\af|dt}{|i\ld+e^{i\rho}|^2}.\EE

From (13) and (14) one readily sees that $\Lim{x\to 0}{\sqrt{\frac{\ld^3}{x}}}I_2=0$ and $\Lim{x\to 0}{\sqrt{\frac{\ld}{x^3}}}I_4=\ld\pi/f(\ld)$. This finally proves that Equations (9) have a simultaneous and positive solution on a curve $y(x,\ld(x),\rho)$ and $\af=\af(y(x))$, for positive $x$ in a neighbourhood of zero.
\\
\\
\centerline{\bf 6. Embeddedness of the Fundamental Piece}
\\

As mentioned in the Introduction, one can profit from the limit-method to prove embeddedness by arguments similar to \cite{V3}, pp489-492. The convergences already studied in Section 5 will be now useful to simplify our task. There we proved the existence of a positive $\eps$ and a curve $\bar{y}:(0,\eps)\to\D^-$, for which the choice $c=c(x,y(x),\af(x))$ in (6) satisfies Equations (9). From now on, $c$ will always represent such a choice.

As a matter of fact, a free parameter $\rho$ can be chosen in $(-\pi/2,0]$, which establishes the curve $\ld:(0,\eps)\to(1,\infty)$, and $y$ is finally given by (10). By taking $z=\frac{x\zeta}{\zeta+i\ld}$, $\zeta\in\B$ and $dH=\ld dh/x$, when $x=0$ the constants $C_{1,2}$ from the previous section will close the periods of the surfaces in \cite{V3}.
 
Choose any $x\in(0,\eps)$ and consider the minimal immersion $X_x:\oM\setminus z^{-1}(\{y^{\,\pm},\bar{y}^{\,\pm}\})\to\R^3$ defined by (5-7). Hence, $X_x$ restricted to $z^{-1}(\D^-)$ can be viewed as a bivalent function $X_x:D^-\setminus\{\bar{y}\}\to\R^3$. Indeed, each branch $w$ of square root in (5) takes any point $q\in D^-\setminus\{\bar{y}\}$ to a pair of points in $\R^3$, say $X_x(q)^+$ and $X_x(q)^-$. Fix $X_x(0)$ as the origin, then one is the image of the other by 180$^\circ$-rotation about $Ox_3$ (see Figure 6a). Moreover, for any closed curve homotopic to $\deh\D^-$, the period vector on it is zero, as proved in Section 5.

Consider the fundamental piece $P$ of $M$. Let $P^-$ be the image of $D^-\setminus\{\bar{y}\}$ in $\R^3$ under $X_x$, and $P^+$ the image in $\R^3$ of $P^-$ under 180$^\circ$-rotation about either $X_x([0,x])$ or $X_x([-1,0])$. Therefore, $P=P^-\cup P^+$. The image of $\deh\D^-$ by $X_x$ is depicted in Figure 6b. 

\begin{figure}[ht]
	\centering
	\includegraphics[scale=0.6]{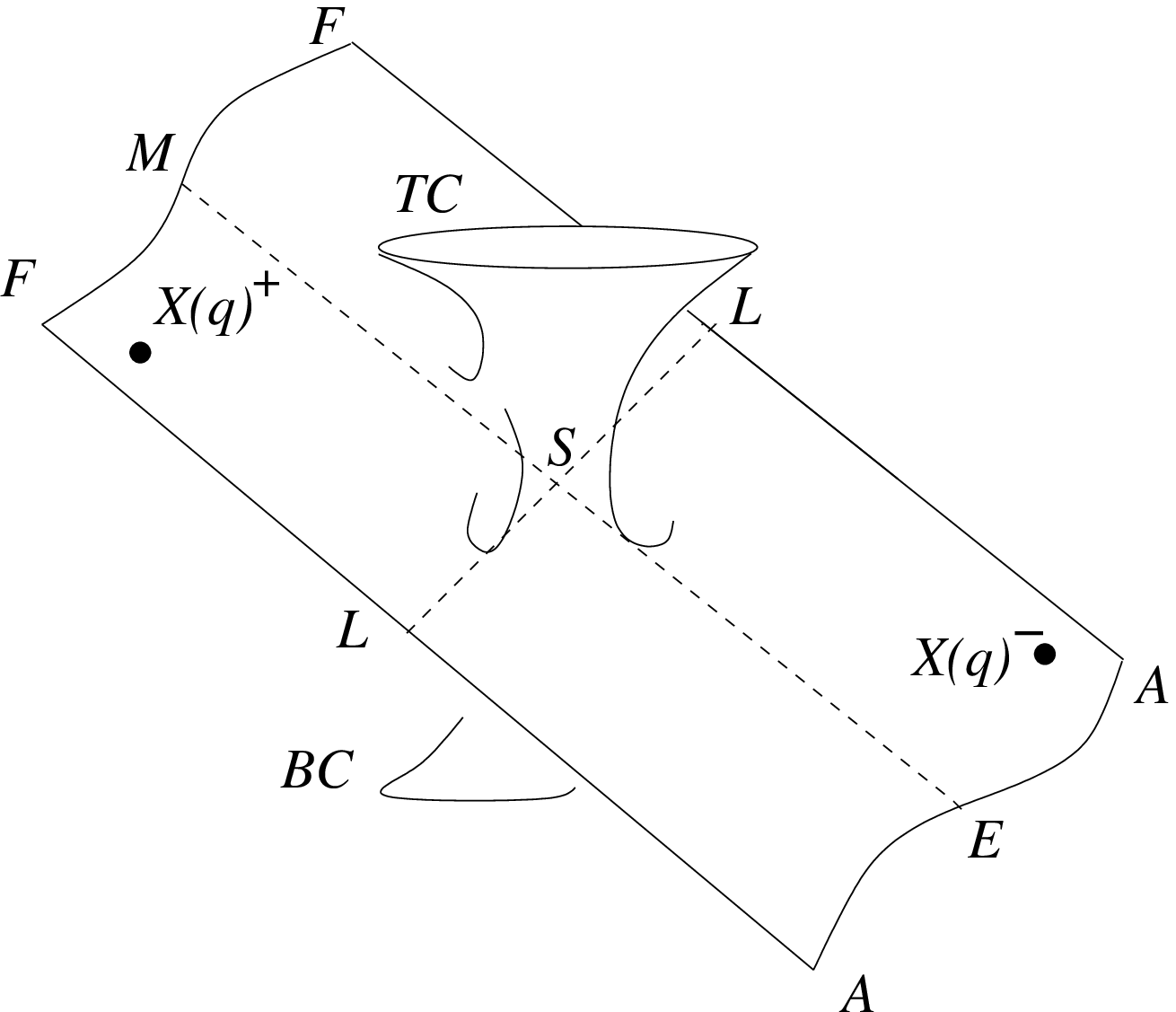}\hspace{2cm} 
        \includegraphics[scale=0.6]{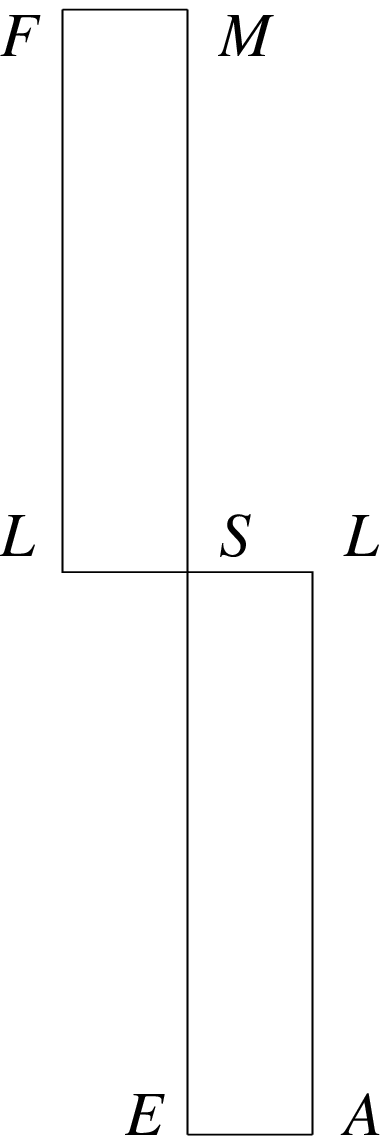}
	\caption{(a) The fundamental piece $P$; (b) the image of $\deh\D^-$ by $X_x$.}
\end{figure}

In Section 5 we defined the set $\B$. Let $\K$ be a compact subset of $\B$ such that $\B\setminus\K=V_{AE}\cup V_{BC}$, where $V_{AE}$ and $V_{BC}$ are connected neighbourhoods of $-i\ld$ and $e^{i\rho}$, respectively. From (11) and (12), one sees that our data $(g,\ld dh/x)$ converge uniformly on $\K$ to the Weierstra\ss  par of the embedded surface $M_{\rm L_b}$, with $\Lim{x\to 0}\sqrt{x\ld}c(x)=c_{\rm L_b}$ (see \cite{V3}). We denote its corresponding minimal embedding by $X_0$.

By \cite{HK} or \cite{MS}, if $V_{BC}$ is small enough, then $X_x\big|_{V_{BC}}$ is the graph of
\[
   x_3(x_1,x_2)=\frac{\eta}{2}\ln(x_1^2+x_2^2)+\mu+\frac{ax_1+bx_2}{x_1^2+x_2^2}+
   \mathcal{O}((x_1^2+x_2^2)^{-1}),
\]
where $\eta,\mu,a$ and $b$ are real numbers. This characterises $X_x\big|_{V_{BC}}$ as a catenoidal end. Since all the parameters vary continuously, $x\to 0$ implies that $X_x\big|_{V_{AE}}$ approaches a pair of Scherk ends, which are also graphs (see \cite{T1}). Thus we can choose $V_{AE}$ small enough for $X_0\big|_{V_{AE}}$ to be a graph. If $x$ is sufficiently close to zero, then the projection of $X_x\big|_{\deh V_{AE}}$ into $x_3=0$ will be a pair of simple closed curves $\mC^{\pm}$, each one consisting of a regular arc and three segments, according to Figure 7. In fact, these curves determine two open regions in the plane, $R^+$ and $R^-$, bounded and simply connected.

\begin{figure}[ht]
	\centering
	\includegraphics[scale=0.8]{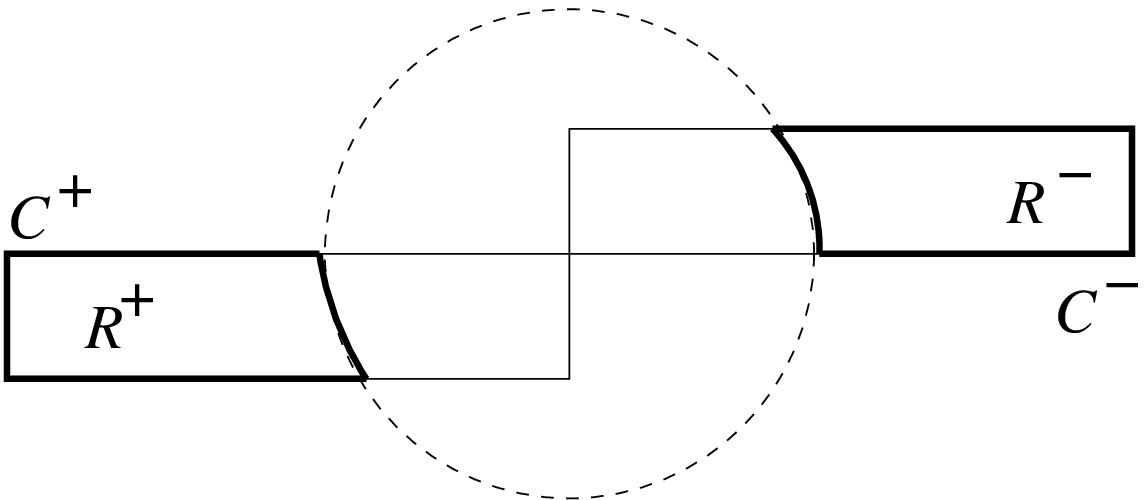}
	\caption{Regions $R^{\pm}$ and curves $\mC^{\pm}$.}
\end{figure}

For sufficiently small $V_{AE}$, $g(V_{AE})$ is contained in a half-sphere of $S^2$. This implies that $(x_1,x_2)\big|_{V_{AE}}$ is an immersion, namely into $R^{\pm}$ because $x_2$ is bounded for any $x\in(0,\eps/2)$. Since $\deh R^{\pm}$ are the monotone curves $\mC^{\pm}$, then $X_x\big|_{\deh V_{AE}}$ is a graph of $x_3$ as a function of $(x_1,x_2)$. The ends do not intersect for sufficiently small $V_{AE}$, $V_{BC}$ and $x$. Hence $X_x\big|_{V_{AE}}$, $X_x\big|_{V_{BC}}$ and $X_x\big|_\K$ are disjoint embeddings in $\R^3$. 

$X_0\big|_\K$ is a compact embedded minimal surface with boundary in $\R^3$. Since the boundary has no self-intersections, then $X_x\big|_\K$ is still embedded for $x$ close enough to zero. Moreover, $X_x\big|_\K$ intersects neither $X_x\big|_{V_{AE}}$ nor $X_x\big|_{V_{BC}}$, otherwise there would be a ball in $\R^3$ containing the whole boundary of $X_x\big|_\K$, but not all the rest of it. This is impossible in the minimal case. Hence, the embedded pieces $X_x\big|_{V_{AE}}$, $X_x\big|_{V_{BC}}$ and $X_x\big|_\K$ make up a whole embedded minimal surface $X_x:\B\to\R^3$, for any $x$ sufficiently close to zero.

This extends to every $x$ in $(0,\eps)$ by means of the maximum principle. Therefore, $P^-$ is embedded in $\R^3$, and since $P^+$ is its image under $180^\circ$ rotation around either $X_x([0,x])$ or $X_x([-1,0])$ of $P^-$, the whole piece $P$ has no self-intersections. The immersion is proper, so $P$ is embedded in $\R^3$.
\\
\\
\\
\centerline{\bf 7. Comments and Remarks}
\\

As we have just seen, the $\CL$-examples build a continuous two-parameter family of minimal surfaces. We claim that, for each example, both parameters control the $Q$-ratio and the logarithmic growth of the ends. From (7) and (10) it is easy to compute
\BE
   Res(dh,z=\bar{y})=\frac{\bar{Y}-2\cos\af}{2(\bar{y}-1/\bar{y})Im\{Y\}}.
\EE

In Section 5 we proved that the fundamental piece has no periods. Consequently, $Re\{2\pi iRes(dh,z=\bar{y})\}=0$ and the following relation holds:
\BE
   \cos\af=\frac{2Re\{y\}}{1+|y|^2}.
\EE
This means that the $y(x)$-curve matches (9) and (16). By substituting (16) in (15) we get
\BE
   Res(dh,z=\bar{y})=\frac{1-|y|^2}{2Im\{Y\}}.
\EE
From (10) and (17) one easily gets $\Lim{x\to 0}{\frac{\ld}{x}Res(dh,z=\bar{y})}=\frac{1}{2}\sec\rho$. This is exactly the value in \cite{V3}, p486. Of course, $x$ and $\rho$ act simultaneously for the logarithmic growth and the $Q$-ratio. In order to get parameters which would control each ratio separately, say $\ell$ and $\boldsymbol{q}$, one should invert the following system of equations:
\[
   \ell=\frac{1-|y(x,\rho)|^2}{2Im\{Y(x,\rho)\}}\eh\eh{\rm and}\eh\eh
   \boldsymbol{q}=\frac{I_1(x,\rho)+c^2(x,\rho)I_3(x,\rho)}{I_6(x,\rho)+c^2(x,\rho)I_7(x,\rho)}.
\]  

Let us now briefly discuss another limit-surface that could be used to obtain the $\CL$-examples. In \cite{V3}, pp492-5, one studies the so-called $M_{\rm C_b}$-surfaces. If we let $x\to 1$, then
\BE
   \lim_{x\to 1}g^2=\frac{-c^2(Z^2-2Re\{Y\}+|Y|^2)}{(Z-2)(Z-2\cos\af)^2}
\EE
and
\BE
   \lim_{x\to 1}dh=\frac{-i(Z-2\cos\af)dz/z}{Z^2-2Re\{Y\}+|Y|^2}.
\EE
From \cite{V3}, p493, one readily recognises (18) and (19) as the Weierstra\ss  data of the surfaces $M_{\rm C_b}$. Moreover, $\Lim{x\to 1}{2I_6}=0$ and $\Lim{x\to 1}{2I_7}=\pi/|2-Y|$, namely $c_{1,2}$ coincide at $x=1$ with the corresponding L\'opez-Ros parameters (41) and (42) of \cite{V3}. Indeed,
\[
   I_6=\int_x^1|dh/w|=\frac{1}{\sqrt{x}}\int_x^1
   \frac{\sqrt{(x^2+1)t-t^2x-x}(t^2+1-2t\cos\af)^2 dt}
   {\sqrt{t}[(t^2+1)^2-2Re\{Y\}(t^3+t)+|Y|^2t^2]^{3/2}}.
\]
By means of the change $v=(x-t)/(x-1)$, it is clear that $\Lim{x\to 1}{2I_6}=0$. 
\[
   2I_7=2\int_x^1\frac{\sqrt{xt}\,dt}
   {\{(1-xt)(t-x)[(t^2+1)^2-2Re\{Y\}(t^3+t)+|Y|^2t^2]\}^{1/2}}.
\]

After applying the change $t=x+(1-x^2)u^2$, it follows that
\[
   2I_7=4\sqrt{x}\int_0^{\frac{1}{\sqrt{1+x}}}
   \frac{\sqrt{(1-x^2)u^2+x}\,du}{\{(1-xu^2)[(t^2+1)^2-2Re\{Y\}(t^3+t)+|Y|^2t^2]\}^{1/2}}.
\]

Therefore,
\[ 
   \lim_{x\to 1}2I_7=4\int_0^{\frac{1}{\sqrt{2}}}
   \frac{du}{[(1-u^2)(4-4Re\{Y\}+|Y|^2)]^{1/2}}=\frac{\pi}{|2-Y|}.
\]

In \cite{V3}, p494, one gets a whole solution curve $(\af(t),y(t))$ on which the equality $c_1|_{x=1}=c_2|_{x=1}$ holds. This curve is obtained by a transversal crossing of graphs on an open subset of $\R^2$. The {\it limit-method} then gives functions $\af(t,x)$ and $y(t,x)$, for $x\in(1-\eps,1)$, which correspond to a simultaneous and positive solution of (9).

\end{document}